\documentclass[12pt, a4paper]{article}

\usepackage{amsmath}
\usepackage{amsopn}
\usepackage{amsfonts}
\usepackage{amsthm}
\usepackage{epsfig, psfrag}

\numberwithin{equation}{section}

\theoremstyle{definition}

\theoremstyle{remark}

\newtheoremstyle{mytheorem}{0.5cm}{0.2cm}{\slshape}{ }{\bfseries}{.}{ }{}
\theoremstyle{mytheorem}
\newtheorem{theorem}{Theorem}[section]

\newtheorem{lemma}[theorem]{Lemma}

\newcommand{\eps}{\varepsilon}

\newcommand{\sN}{{\mathcal{N}}}
\renewcommand{\SS}{{\mathbb{S}}}

\newcommand{\ti}{{\to \infty}}
\newcommand{\Prob}[1]{\mathbf{P}\{#1\}}
\newcommand{\E}{\mathbf{E}}
\newcommand{\Ind}{\mathbf{1}}
\newcommand{\standardspace}[1]{\mathbb{#1}}
\newcommand{\R}{\standardspace{R}}
\newcommand{\Sphere}{\standardspace{S}^{d-1}}
\newcommand{\Beta}{\mathrm{B}}
\newcommand{\rhoH}{\rho_{\mathrm{H}}}
\newcommand{\thf}{\frac{1}{2}}
\newcommand{\Pit}{\tilde{\Pi}}
\newcommand{\xit}{\tilde{\xi}}
\DeclareMathOperator{\diam}{diam}
\DeclareMathOperator{\supp}{supp}

\newcommand{\as}{\quad\text{as}\;}

\usepackage{fancybox}
\newcommand{\query}[1]{\medskip \noindent
  \shadowbox{\begin{minipage}[t]{\textwidth} {#1}
    \end{minipage}}\medskip}
\renewcommand{\query}[1]{}

\begin{document}

\bibliographystyle{plain}

\title{A limit theorem for the maximal interpoint distance of a random
  sample in the unit ball\footnote{Supported by Swiss National
    Foundation Grant No.  200021-103579}}

\author{Michael Mayer, Ilya Molchanov \\
  \textit{\normalsize Department of Mathematical Statistics
    and Actuarial Science},\\
  \textit{\normalsize University of Bern, Sidlerstrasse 5, CH-3012
    Bern, Switzerland}}

\date{}

\maketitle

\begin{abstract}
  We prove a limit theorem for the the maximal interpoint distance
  (also called the diameter) for a sample of $n$ i.i.d. points 
  in the unit $d$-dimensional ball for $d\geq 2$. The exact form of
  the limit distribution and the required normalisation are derived
  using assumptions on the tail of the interpoint distance for two
  i.i.d.  points. The results are specialised for the cases when the
  points have spherical symmetric distributions, in particular, are
  uniformly distributed in the whole ball and on its boundary.


  \noindent AMS 60D05; 60G55; 60G70
\end{abstract}

\section{Introduction}
\label{sec:introduction}

Asymptotic behaviour of random polytopes formed by taking convex hulls
of samples of i.i.d. points has been thoroughly investigated in the
literature, see, e.g., \cite{gru93a,schn88ap} for surveys of classical
results and \cite{reit05} for more recent studies. Consider a
\emph{random polytope} $P_n$ obtained as the convex hull of $n$ i.i.d.
points $\xi_1,\dots,\xi_n$ sampled from the Euclidean space $\R^d$.

Most of results about random convex hulls are available in the planar
case, i.e. for $d=2$.  The typical questions about random polytopes
$P_n$ concern the limit theorems for the geometric characteristics of
$P_n$, for instance the area, the perimeter and the number of vertices
of $P_n$, see \cite{bra:hsi98, gro88,schn88ap}. Further important
results concern the quantities that characterise the \emph{worst case}
approximation, notably the Hausdorff distance between $K$ and $P_n$,
see \cite{bra:hsi:bin98,dum:wal96}. It is well known \cite{schn} that
the Hausdorff distance between two convex sets equals the uniform
distance between their support functions defined on the unit sphere,
i.e.
\begin{displaymath}
  \rhoH(P_n,K)=\sup_{u:\;\|u\|=1}(h(K,u)-h(P_n,u))\,,
\end{displaymath}
where $\|u\|$ is the Euclidean norm of $u\in\R^d$,
\begin{displaymath}
  h(K,u)=\sup\{\langle u,x\rangle:\; x\in K\}
\end{displaymath}
is the support function of $K$ (and similar for $P_n$) and $\langle
u,x\rangle$ is the scalar product in $\R^d$. For instance
\cite{dum:wal96} shows that for uniformly distributed points
$\rhoH(P_n,K)$ is of order ${\mathcal{O}}((n^{-1}\log n)^{2/(d-1)})$
if $K$ is sufficiently smooth.

The results on the \emph{best case} approximation concern the
behaviour of the infimum of the difference between $h(K,u)$ and
$h(P_n,u)$.  One of the few results in this direction states that if
$K$ is smooth, then $n(h(K,u)-h(P_n,u))$ (as a stochastic process
indexed by $u$ from the unit sphere $\Sphere$) epi-converges in
distribution to a certain process derived from the Poisson point
process on $\Sphere\times[0,\infty)$, see \cite{mo95} and
\cite[Th.~5.3.34]{mo1}.  The epi-convergence implies the weak
convergence of infima on each compact set. In particular,
$n\inf_{u\in\Sphere} (h(K,u)-h(P_n,u))$ converges in distribution to
an exponentially distributed random variable, i.e. the best
approximation error is of the order of $n^{-1}$. If the points are
uniformly distributed in $K$, then this exponential random variable
has the mean being the ratio of the volume of $K$ and its surface
area, see \cite[Ex.~5.3.35]{mo1}.  Further results along these lines
can be found in \cite{schr02a}.

The best case approximation can be also studied by considering how
fast the diameter of $P_n$, $\diam P_n$, approximates $\diam K$. By
\emph{diameter} we understand the maximum distance between any two
points from the set. Note that $\diam K$ is not necessarily equal to
the diameter of the smallest ball that contains $K$. This is the case,
e.g. if $K$ is a triangle.

A limit theorem for the diameter of $P_n$ was proved in
\cite{app:naj:rus02} for uniformly distributed points in a compact set
$K$ with unique longest chord (whose length is the diameter) and such
that the boundary of $K$ near the endpoints of this major axis is
locally defined by regularly varying functions with indices strictly
larger than $0.5$. These assumptions are fairly restrictive and
exclude a number of interesting smooth sets $K$, in particular balls
and ellipsoids.

For $K$ being the unit disk on the plane, \cite{app:naj:rus02}
provides only bounds for the limit distribution, even without proving
the existence of the limit. In particular, \cite[Th.~4]{app:naj:rus02}
states that
\begin{align}
  \label{eq:ap-r}
  1-\exp\left\{-\frac{4t^{5/2}}{3^{5/2}\pi}\right\}
  &\le \liminf_{n \ti} \Prob{n^{4/5}(2-\diam P_n) \le t} \notag \\
  &\le \limsup_{n \ti} \Prob{n^{4/5}(2-\diam P_n) \le t} \notag \\
  &\le 1 - \exp\left\{-\frac{4t^{5/2}}{\pi}\right\}\,, \quad t>0\,.
\end{align}

In the classical theory of extreme values it is possible either to
normalise the maximum of a random sample by dividing or multiplying
its (possibly shifted or translated) maximum with normalising
constants that grow to infinity. The first case corresponds to
samples with possibly unbounded values, while the second one appears
if samples with a finite right end-point of the distribution are
considered. Quite similarly, in the extreme problems for random
polytopes one can consider samples supported by the whole $\R^d$ or by
a compact convex subset $K$ in $\R^d$. In this paper we consider only
the latter case.  The limit theorems for the largest interpoint
distances for samples from the whole $\R^d$ have been proved in
\cite{mat:ruk93} for the normally distributed samples and in
\cite{hen:klein96} for more general spherically symmetric samples.

In this paper we state limit laws of the diameters of $P_n$, where
$P_n$ is the convex hull of a sample $\Xi_n = \{\xi_1,\dots,\xi_n\}$
of independent points distributed in the $d$-dimensional unit ball
\begin{displaymath}
  B=\{x\in \R^d: \|x\| \le 1\}
\end{displaymath}
according to some probability measure $\kappa$. The diameter of a set
$F\subset \R^d$ is determined by its largest interpoint distance, i.e.
by
\begin{displaymath}
  \diam(F) = \sup_{x,y \in F}\|x-y\|\,,
\end{displaymath}
and it is obvious that the diameter of $F$ equals the diameter of its
convex hull. In the special case when $\kappa$ is the uniform
distribution, the following result provides a considerable improvement
of \cite[Th.~4]{app:naj:rus02}.

\begin{theorem}
  \label{theoremuniform}
  As $n\ti$, the diameter of the convex hull $P_n$ of $n$
  independent points distributed uniformly on the $d$-dimensional unit
  ball $B$, $d\ge 2$, has limit distribution given by
  \begin{multline*}
    \Prob{n^{\frac{4}{d+3}}(2 - \diam P_n) \le t } \to 1 -
    \exp\left\{-\frac{2^d d \Gamma(\frac{d}{2}+1)}
      {\sqrt{\pi}(d+1)(d+3)\Gamma(\frac{d+1}{2})}
      t^{\frac{d+3}{2}}\right\},\\ \quad t > 0,
  \end{multline*}
  where $\Gamma(x)=\int_0^\infty s^{x-1}e^{-s}ds$ denotes the Gamma
  function.
\end{theorem}

This theorem is proved by showing that the same limit distribution is
shared by the diameter of a homogeneous Poisson process $\Pi$ of
constant intensity $\lambda = n/\mu_d(B)$ restricted on $B$, so that
the total number of points in $\Pi$ has mean $n$. See
Section~\ref{sec:de-poissonisation} for a more general
de-Poissonisation argument, which implies that the diameter of a
general binomial process with $n$ points and of the corresponding
Poisson process share the same limiting distributions (if the limit
distribution exists). 

The problem in dimension 1 is very easy to solve, see e.g.
\cite{gal78}. It is interesting to note that if all ${n \choose 2}$
random distances $\|\xi_i-\xi_j\|$ are treated as an i.i.d.  sequence
with the common distribution determined by the length of the random
chord in $K$, then the maximum of these distances has the same limit
law as described in Theorem~\ref{theoremuniform}.  This is explained
by the fact that only different pairs of points contribute to $\diam
P_n$, while the probability that a point has considerably large
interpoint distances with two or more other points is negligible. This
argument stems from \cite{sil:bro78} and was used in the proofs in
\cite{mat:ruk93} and \cite{hen:klein96}. Our proof relies on
properties of the Poisson process with a subsequent application of a
de-Poissonisation argument.

In Section~\ref{sec:diam-poiss-proc} we establish the asymptotic
behaviour of the diameter for a Poisson point process in $B$ with
growing intensity. The conditions on the intensity $\kappa$ of the
Poisson point process require certain asymptotic behaviour of the
distance between two typical points of the process and a certain
boundedness condition on $\kappa$. For instance, these conditions are
fulfilled in the uniform case.

In Section~\ref{sec:spherically-symmetric} we investigate the
asymptotic behaviour of the diameter of $\Pi_{n\kappa}$, where
$\kappa$ is a spherically symmetric distribution.
Section~\ref{sec:examples} describes several examples, in
particularly, where $\kappa$ is the uniform measure on $B$ and on
$\Sphere$, respectively.  Further examples concern distributions which
are not spherically symmetric.

\bigskip

The ball of radius $r$ centered at the origin is denoted by $B_r$.  By
$\mu_d$ we denote the $d$-dimensional Lebesgue measure in $\R^d$.
Furthermore, $\mu_{d-1}$ is the surface area measure on the unit
sphere $\Sphere$. By $\kappa$ we understand a certain fixed
probability measure on $B$ and $\xi_1,\xi_2,\ldots$ are i.i.d. points
distributed according to $\kappa$.

For any set $F$ in $\R^d$, $\check{F}$ denotes the reflected set
$\{-x:\; x\in F\}$ and $\tilde{F}$ is the corresponding difference set
\begin{displaymath}
  \tilde{F}=F+\check{F}=\{x-y:\; x,y\in F\}\,.
\end{displaymath}

Finally, the letter $\Pi_\nu$ stands for the Poisson process on $B$ of
intensity measure $\nu$, where we write shortly $\Pi$ if no ambiguity
occurs or the intensity measure is immaterial. Note that $\Pi(F)$
denotes the number of points of a point process inside a set $F$, so
that $\Pi(F)=0$ is equivalent to $\Pi\cap F=\emptyset$.

\section{Diameters for Poisson processes}
\label{sec:diam-poiss-proc}

Consider a Poisson process $\Pi=\Pi_{n\kappa}$ in the unit ball $B$
with the intensity measure proportional to a probability measure
$\kappa$ on $B$. Consider the convolution of $\kappa$ with the
reflected $\kappa$, i.e. the probability measure $\kappa$ that
determines the distribution of $\xit=\xi_1-\xi_2$ for i.i.d.
$\xi_1,\xi_2$ distributed according to $\kappa$. Assume throughout
that the support of $\tilde\kappa$ contains points with norms
arbitrarily close to $2$. In this case the diameter of $\Pi_{n\kappa}$
approaches $2$ as $n\ti$. In this section we determine the asymptotic
distribution of $2-\diam(\Pi_{n\kappa})$ as $n\ti$.

The distribution of the diameter of $\Pi$ is closely related to the
probability that the inner $s$-shell $B_2\setminus B_{2-s}$ of the
ball of radius $2$ contains no points of $\Pit=\Pi+\check{\Pi}$.
Indeed
\begin{displaymath}
  \Prob{\diam \Pi \le 2-s} = \Prob{\Pit(B_2\setminus
  B_{2-s})=0}\,,
\end{displaymath}
and by the symmetry of $\Pit$,
\begin{equation}
  \label{eq:hsp}
  \Prob{\diam \Pi \le 2-s}
  = \Prob{\Pit((B_2\setminus B_{2-s})\cap H)=0}\,,
\end{equation}
where $H$ is any halfspace bounded by a $(d-1)$-dimensional hyperplane
passing through the origin.

For each $A\subset\Sphere$ define
\begin{equation}
  \label{eq:a-def}
  A_s=\{r x:\; x \in A, \; r\in[2-s,2]\}\,.
\end{equation}
For each point $u\in\Sphere$ define a cap of the unit ball of height
$s\in(0,1)$ by
\begin{displaymath}
  D_s(u) = B\cap \{x \in \R^d:\; \langle x, u\rangle \ge 1-s\}\,,
\end{displaymath}
where $\langle x, u\rangle$ denotes the scalar product.  For
$A\subset\Sphere$ define
\begin{displaymath}
  D_s(A) = \cup_{u \in A}D_s(u)\,, \quad s\in(0,1)\,.
\end{displaymath}
Then $D_s(A)$ and $D_s(\check{A})$ are subsets of $B\setminus B_{1-s}$
with the property that $x_1-x_2\in A_s$ for some $x_1,x_2\in B$
implies that $x_1$ belongs to $D_s(A)$ and $x_2$ to $D_s(\check{A})$.

\begin{lemma}
  \label{le:candidate-set}
  For each $A\subset\Sphere$, $s\in(0,1)$ and each $x_1\in B\setminus
  D_s(A)$ and $x_2\in B$, we have $x_1-x_2\notin A_s$.
\end{lemma}
\begin{proof}
  By definition of $D_s(u)$ and the fact that $\|x_2\|\le 1$, the
  inequality
  \begin{displaymath}
    \langle u,x_1-x_2\rangle = \langle u,
    x_1\rangle + \langle u,-x_2 \rangle < 2 - s
  \end{displaymath}
  holds for each $u\in A$.  If $x_1-x_2\in A_s$, then $\|x_1-x_2\|\ge
  2-s$ and $u_0 = (x_1-x_2)\|x_1-x_2\|^{-1}\in A$. Now write $2-s \le
  \|x_1-x_2\| = \langle u_0,x_1-x_2 \rangle$, which is a contradiction
  to the first inequality, and hence the claim.
\end{proof}

\begin{lemma}
  \label{lem:dsb}
  For each $s\in(0,1)$ and $A\subset\Sphere$, the set $D_s(A)$ lies
  inside the $\sqrt{2s}$-neighbourhood of $A$.
\end{lemma}
\begin{proof}
  Consider arbitrary $u\in A$. Since
  \begin{displaymath}
    \|x-u\|^2 = \|x\|^2 + \|u\|^2-2\langle x, u \rangle
    \leq 2 - 2(1 - s) = 2s\,,
  \end{displaymath}
  every point $x\in D_s(u)$ is located within distance at most
  $\sqrt{2s}$ from $u$.
\end{proof}

The following lemma follows from Lemma~\ref{le:candidate-set} and the
independence property of the Poisson process.

\begin{lemma}
  \label{le:independence}
  For any $A\subset\Sphere$ and $s\in(0,1)$,
  \begin{equation}
    \label{eq:pie}
    \Prob{\Pit \cap A_s \ne \emptyset}
    = \Prob{\Pit \cap A_s \ne \emptyset, \Pi \cap D_s(A) \ne \emptyset,
      \check{\Pi} \cap D_s(A) \ne\emptyset}\,.
  \end{equation}
  If $A',A''\subset\Sphere$ and
  \begin{displaymath}
    D_s(A')\cap D_s(A'') = D_s(A')\cap D_s(\check{A''})=\emptyset\,,
  \end{displaymath}
  then the random variables $\Pit(A'_s)$ and $\Pit(A''_s)$ are
  independent.
\end{lemma}

The following lemma bounds $\Prob{\Pit_{n\kappa}\cap A_s\ne\emptyset}$
using $\Prob{\xi_1-\xi_2\in A_s}$ for independent points $\xi_1$ and
$\xi_2$ distributed according to $\kappa$.

\begin{lemma}
  \label{le:large-block}
  For each $A\subset\Sphere$ and $0<s<1$, we have
  \begin{align*}
    n^2e^{-n(a+\check{a})}\Prob{\xit\in A_s}
    & \le \Prob{\Pit_{n\kappa}\cap A_s \ne \emptyset} \\
    &\le n^2(1+na\check{a}(a+\check{a})) \Prob{\xit\in A_s}\,,
  \end{align*}
  where $a=\kappa(D_s(A))$, $\check{a}=\kappa(D_s(\check{A}))$ and
  $\xit=\xi_1-\xi_2$ for $\xi_1$ and $\xi_2$ being independent
  points distributed according to $\kappa$.
\end{lemma}
\begin{proof}
  Let $\zeta_1$ and $\zeta_2$ be Poisson distributed with means $na$
  and $n\check{a}$ respectively, so that $\zeta_1$ and $\zeta_2$
  represent the numbers of points of $\Pi$ in $D_s(A)$ and
  $D_s(\check{A})$ respectively.  First, (\ref{eq:pie}) implies that
  \begin{align*}
    \Prob{\Pit \cap A_s \ne \emptyset}
    &= \Prob{\Pit(A_s) \geq 1, \zeta_1 \geq 1,\zeta_2 \geq 1}\\
    &\ge \Prob{\Pit(A_s) = 1, \zeta_1=1,\zeta_2=1}\,.
  \end{align*}
  An upper bound follows from
  \begin{align*}
    \Prob{\Pit \cap A_s \ne \emptyset}
    &= \Prob{\Pit(A_s) \geq 1, \zeta_1 \geq 1,\zeta_2 \geq 1 }\\
    &\leq \Prob{\Pit(A_s) = 1, \zeta_1=1,\zeta_2=1}+I\,,
  \end{align*}
  where
  \begin{align*}
    I&=\sum_{k_1,k_2=2,\atop \max(k_1,k_2)\geq2}^\infty
    \Prob{\Pit(A_s)\geq1,\zeta_1=k_1,\zeta_2=k_2}\,.
  \end{align*}
  The subadditivity of probability and the fact that $\zeta_1$ and
  $\zeta_2$ are independent immediately imply that
  \begin{displaymath}
    \Prob{\Pit(A_s)\geq1|\zeta_1=k_1,\zeta_2=k_2}
    \leq k_1k_2\Prob{\xi_1-\xi_2\in A_s}\,.
  \end{displaymath}
  Thus,
  \begin{align*}
    I&\leq \Prob{\xi_1-\xi_2\in A_s}
    \left(\E(\zeta_1\zeta_2)-\Prob{\zeta_1=1}\Prob{\zeta_2=1}\right)\\
    &=\Prob{\xi_1-\xi_2\in
      A_s}(n^2a\check{a}-n^2a\check{a}e^{-n(a+\check{a})})\\
    &\leq \Prob{\xi_1-\xi_2\in A_s}n^3a\check{a}(a+\check{a})\,.
  \end{align*}
  Now write
  \begin{align*}
    \Prob{\Pit(A_s) = 1, \zeta_1=1,\zeta_2=1} &= \Prob{\Pit(A_s) = 1\
      | \zeta_1=1,\zeta_2=1}
    n^2a\check{a}e^{-n(a+\check{a})}\\
    &= \Prob{\eta_1-\eta_2 \in A_s} n^2a\check{a}e^{-n(a+\check{a})}\,,
  \end{align*}
  where $\eta_1$ and $\eta_2$ are independent points distributed
  according to the normalised measure $\kappa$ restricted onto
  $D_s(A)$ and $D_s(\check{A})$ respectively. Because of
  Lemma~\ref{le:candidate-set},
  \begin{displaymath}
    \Prob{\eta_1-\eta_2 \in A_s}=\frac{1}{a\check{a}}
    \Prob{\xi_1-\xi_2\in A_s}\,,
  \end{displaymath}
  and the proof is complete.
\end{proof}

Let
\begin{displaymath}
  C(u,r) = \{x\in\Sphere: \|x-u\| \leq r\}\,,\quad u\in\Sphere,\; r>0\,,
\end{displaymath}
denote the spherical ball, so that
\begin{displaymath}
  C_s(u,r) = \{rx: x \in C(u,r), r \in [2-s,2]\}
\end{displaymath}
in accordance with (\ref{eq:a-def}).

Introduce the following assumption on the distribution of the
difference $\xit$ between two independent points in $B$ distributed
according to $\kappa$. Assume that for a finite non-trivial measure
$\sigma$ on $\Sphere$, some $\gamma>0$ and $[\delta',\delta'']\subset
(0,\thf)$ we have
\begin{equation}
  \label{eq:su}
  \lim_{s\downarrow0} \frac{\Prob{\xit\in C_s(u,z_s)}}
  {s^\gamma\sigma(C(u,z_s))} =1 
\end{equation}
and 
\begin{equation}
  \label{eq:ba}
    \lim_{s\downarrow 0} s^{-\gamma/2} \kappa(D_s(C(u,z_s))) = 0
\end{equation}
uniformly in $u\in\Sphere$ and $z_s\in [s^{\delta'},s^{\delta''}]$. If
$u$ does not belong to the support of $\sigma$, then the denominator
in (\ref{eq:su}) equals zero for all sufficiently small $s$, and
(\ref{eq:su}) is understood as the fact that the numerator also equals
zero for all sufficiently small $s$.  Since $\xit$ has a centrally
symmetric distribution, the measure $\sigma$ is necessarily centrally
symmetric.

\begin{lemma}
  \label{lem:b-dens}
  If (\ref{eq:su}) holds with $\gamma<d+1$, $\kappa$ is absolutely
  continuous on $B_1\setminus B_{1-s}$ for some $s>0$ and possesses
  there a bounded density, then (\ref{eq:ba}) holds with 
  \begin{equation}
    \label{eq:gamma-delta}
    \frac{\gamma-2}{2(d-1)}< \delta'\leq \delta''<\thf\,.
  \end{equation}
\end{lemma}
\begin{proof}
  It suffices to show that, for any given $u\in\Sphere$,
  \begin{displaymath}
    \lim_{s\downarrow 0} s^{-\gamma/2} \mu_d(D_s(C(u,s^\delta))) = 0\,.
  \end{displaymath}
  By Lemma~\ref{lem:dsb}, noticing that $\delta<\thf$, this would
  follow from
  \begin{displaymath}
    s^{-\gamma/2} \mu_{d-1}(C(u,2s^\delta))s\to 0\as s\downarrow0\,.
  \end{displaymath}
  The latter is the case, since $-\thf\gamma+\delta(d-1)+1>0$ for all
  $\delta\in[\delta',\delta'']$. Finally, $\gamma<d+1$ implies that 
  $\frac{\gamma-2}{2(d-1)}<\thf$, so that (\ref{eq:gamma-delta}) makes
  sense. 
\end{proof}

In general, (\ref{eq:ba}) is weaker than the boundedness of the
density of $\kappa$ with respect to the Lebesgue measure, which would,
e.g., exclude the case where $\kappa$ is supported by $\Sphere$.

\begin{lemma}
  \label{lemma:negligible}
  If (\ref{eq:su}) and (\ref{eq:ba}) hold, then, for every measurable
  set $A\subset\Sphere$ and $c>0$,
  \begin{displaymath}
    \lim_{s\downarrow0} \Prob{\Pit_{n\kappa}\cap A_s\neq\emptyset}
    \leq c^2\sigma(A)\,,
  \end{displaymath}
  where $n=cs^{-\gamma/2}$.
\end{lemma}
\begin{proof}
  Cover $A$ by spherical balls $C(u_i,s^\delta)$, $i=1,\dots,m$, of
  diameter $s^\delta$, where $\delta\in[\delta',\delta'']$.  Then
  \begin{displaymath}
    \Prob{\Pit_{n\kappa}\cap A_s\neq\emptyset}
    \leq \sum_{i=1}^m
    \Prob{\Pit_{n\kappa}\cap C_s(u_i,s^\delta)\neq\emptyset}\,.
  \end{displaymath}
  By the choice of $n$, Lemma~\ref{le:large-block} and (\ref{eq:su}),
  \begin{align*}
    \Prob{\Pit_{n\kappa}\cap C_s(u_i,s^\delta)\neq\emptyset} &\leq c^2
    s^{-\gamma}\Prob{\xit\in C_s(u_i,s^\delta)}(1+na_i\check{a}_i)\\
    &\leq c^2(1+\eps)\sigma(C(u_i,s^\delta))(1+na_i\check{a}_i)
  \end{align*}
  for any $\eps>0$ and all sufficiently small $s$, where
  $a_i=\kappa(D_s(C(u_i,s^\delta)))$ and
  $\check{a}_i=\kappa(D_s(C(-u_i,s^\delta)))$.  Condition
  (\ref{eq:ba}) implies that $na_i\check{a}_i\to 0$ as $s\downarrow0$.
  Therefore,
  \begin{displaymath}
    \Prob{\Pit_{n\kappa}\cap C_s(u_i,s^\delta)\neq\emptyset} 
    \leq c^2(1+\eps)^2\sigma(C(u_i,s^\delta))
  \end{displaymath}
  for all sufficiently small $s$. Thus,
  \begin{displaymath}
    \Prob{\Pit_{n\kappa}\cap A_s\neq\emptyset}
    \leq c^2(1+\eps)^2\sum_{i=1}^m \sigma(C(u_i,s^\delta))
  \end{displaymath}
  for all sufficiently small $s$. The statement is proven by taking
  infimum in the right-hand side over all possible ball-coverings of
  $A$.
\end{proof}

In the following we need the following assumption on $\sigma$:
\begin{description}
\item[(S)] $\sigma$ is a measure on $\Sphere$ with finite total mass
  $\sigma_0$ such that 
  \begin{equation}
    \label{eq:sigma-l}
    \sigma(A)\leq f(\mu_{d-1}(A))
  \end{equation}
  for all measurable $A\subset\Sphere$ with a function $f$ such that
  $f(x)\to 0$ as $x\downarrow 0$. 
\end{description}

It is easy to see that (\ref{eq:sigma-l}) holds if $\sigma$ is
absolutely continuous with respect to $\mu_{d-1}$ and has a bounded
density. An atomic $\sigma$ clearly violates \textbf{(S)}.

\begin{theorem}
  \label{thr:non-unif}
  Assume that~(\ref{eq:su}) and (\ref{eq:ba}) hold with
  $\delta'<\delta''$ and a $\sigma$ that satisfies \textbf{(S)}.  Then
  \begin{equation}
    \label{eq:sgamma}
    \lim_{n\ti} \Prob{n^{2/\gamma}(2-\diam(\Pi_{n\kappa}))\leq t}
    =1-e^{-\thf t^\gamma \sigma_0}\,,\quad t\geq0\,,
  \end{equation}
  where $\sigma_0=\sigma(\Sphere)$.
\end{theorem}
\begin{proof}
  Let $\Sphere_+$ denote the half-sphere, obtained by intersection
  $\Sphere$ with any fixed half-space $H$, for instance given by
  (\ref{eq:hsp}).  Fix any $\eps>0$ and consider disjoint spherical
  balls $C(x_i,s^{\delta_i})$, $i=1,\dots,m$, where $x_i\in\Sphere_+$
  and $\delta_i\in[\delta',\delta'']$.
  Since these spherical balls are constructed using varying scales of
  $s$, it is possible to pack them arbitrarily dense as
  $s\downarrow0$, i.e. assume that the Lebesgue measure of the
  uncovered part is smaller than $\eps$. 

  Define the spherical balls
  \begin{displaymath}
    A^{(i)}=C(x_i,s^{\delta_i}-\sqrt{2s})\,, \quad i=1,\dots,m\,.
  \end{displaymath}
  Since $\sqrt{2s}\leq s^{\delta_i}$ for all sufficiently small $s$,
  Lemma~\ref{lem:dsb} implies that $D_s(A^{(i)})$, $i=1,\dots,m$, are
  pairwise disjoint for all sufficiently small $s$. By
  Lemma~\ref{le:independence}, the random variables $\Pit(A_s^{(i)})$,
  $i=1,\dots,m$, are independent.

  Denote
  \begin{equation}
    \label{eq:delta-s}
    \Delta(s)=\Sphere_+\setminus (A^{(1)}\cup\cdots\cup A^{(m)})
  \end{equation}
  to be the uncovered part of $\Sphere_+$ left by the $A^{(i)}$'s.
  The Lebesgue measure of $\Delta(s)$ equals the sum of the
  $\mu_{d-1}$-measure of the part left uncovered by
  $C(x_i,s^{\delta_i})$, $i=1,\dots,m$, and the sum of the measures of
  $C(x_i,s^{\delta_i})\setminus A^{(i)}$. Thus
  \begin{align*}
    \mu_{d-1}(\Delta(s))&\leq \eps + \sum_{i=1}^m
    c_1(d-2)s^{\delta_i(d-2)}q\sqrt{2s}\\
    &\leq \eps +
    \sum_{i=1}^m c_1 s^{\delta_i(d-1)} \sqrt{2s}s^{-\delta'}\\
    &\leq \eps +
    c_2\sqrt{2s}s^{-\delta'}\leq 2\eps
  \end{align*}
  for all sufficiently small $s$, where $c_1$ and $c_2$ are positive
  constants. Since the chosen points $x_1,\dots,x_m$ do not include at
  most $\eps$ of the atomic part of $\sigma$, condition \textbf{(S)}
  implies that $\sigma(\Delta(s))$ is smaller than $\eps+f(2\eps)$ for
  all sufficiently small $s$. In turn, $\eps+f(\eps)$ can be made
  smaller than any given $\eps'>0$. By Lemma~\ref{lemma:negligible},
  \begin{equation}
    \label{eq:d-small}
    \lim_{s\downarrow0}\Prob{\Pit\cap\Delta_s(s)\neq\emptyset}
    \leq t^2\eps'\,.
  \end{equation}
  
  For any fixed $t>0$ consider the Poisson process $\Pi$ with
  intensity measure $n\kappa$ with $n=ts^{-\gamma/2}$ for a fixed $t$.
  Then
  \begin{displaymath}
    \Prob{\diam(\Pi)\leq 2-s}
    =\Prob{\Pit\cap A_s^{(i)}=\emptyset,\, i=1,\dots,m,\;
      \Pit\cap\Delta_s(s)=\emptyset}\,.
  \end{displaymath}
  By the independence of $\Pit(A_s^{(i)})$, $i=1,\dots,m$,
  \begin{displaymath}
    I\leq \Prob{2-\diam(\Pi)\leq s}
    \leq I+\Prob{\Pit\cap\Delta_s(s)\neq\emptyset}\,,
  \end{displaymath}
  where
  \begin{displaymath}
    I=1-\prod_{i=1}^m \Prob{\Pit(A_s^{(i)})=0}\,.
  \end{displaymath}
  By Lemma~\ref{le:large-block},
  \begin{align*}
    \prod_{i=1}^m (1-n^2(1+y_1(s))\Prob{\xit\in A_s^{(i)}})
    &\leq \prod_{i=1}^m \Prob{\Pit(A_s^{(i)})=0}\\
    &\leq \prod_{i=1}^m (1-n^2e^{-2y_2(s)}
    \Prob{\xit\in A_s^{(i)}})\,,
  \end{align*}
  where
  \begin{align*}
    y_1(s)&=\max_{i=1,\dots,m}
    n\kappa(D_s(A^{(i)}))\kappa(D_s(\check{A}^{(i)}))\,,\\
    y_2(s)&=\max_{i=1,\dots,m} n\kappa(D_s(A^{(i)}))\,.
  \end{align*}
  By (\ref{eq:ba}), $y_2(s)$ (and thereupon also $y_1(s)$) converge to
  zero as $s\downarrow 0$ for $n=ts^{-\gamma/2}$.    
  By (\ref{eq:su}) with $z_s=s^{\delta_i}-\sqrt{2s}$,
  \begin{displaymath}
    \frac{n^2\Prob{\xit\in A_s^{(i)}}}{\sigma(A^{(i)})}\to
    t^\gamma
    \as s\downarrow0\,.
  \end{displaymath}
  Since $y_1(s)\to0$ and $y_2(s)\to0$,
  \begin{displaymath}
    \lim_{s\downarrow0} \prod_{i=1}^m (1-n^2(1+y_1(s))
    \Prob{\xit\in A_s^{(i)}})
    =\lim_{s\downarrow0} \prod_{i=1}^m (1-t^\gamma\sigma(A^{(i)}))\,,
  \end{displaymath}
  and
  \begin{displaymath}
    \lim_{s\downarrow0}\prod_{i=1}^m (1-n^2e^{-2y_2(s)}\Prob{\xit\in
    A_s^{(i)}})=\lim_{s\downarrow0} \prod_{i=1}^m
    (1-t^\gamma\sigma(A^{(i)}))\,.
  \end{displaymath}
  By taking logarithm, and using the inequality $|\log(1+x)-x|\leq
  |x|^2$ for $|x|<1$, we see that
  \begin{align}
    \label{eq:inf-product}
    \lim_{s\downarrow 0} \prod_{i=1}^m (1-t^\gamma\sigma(A^{(i)}))
    &=\exp\left\{ -t^\gamma\lim_{s\downarrow0}
      \sum_{i=1}^m \sigma(A^{(i)})\right\}\\
    &=\exp\{-\thf t^\gamma\sigma_0\}\,.\notag
  \end{align}
  For this, note that $\sigma$ is necessarily symmetric, so that
  $\sigma(\Sphere_+)=\sigma_0/2$.  Finally, (\ref{eq:sgamma}) is
  obtained from (\ref{eq:d-small}) and the choice of
  $n=(t/s)^{-\gamma/2}$.
\end{proof}

Instead of imposing \textbf{(S)} it is possible to request that for
every $s>0$ there exists a covering of $\Sphere$ by spherical balls
$C(x_i,s^{\delta_i})$ of radii $s^{\delta_i}$ with
$\delta_i\in[\delta',\delta'']\subset(0,\thf)$ such that
$\sigma(\Delta(s))\to 0$ as $s\downarrow0$, where $\Delta(s)$ is given
by (\ref{eq:delta-s}). Since this condition always holds in dimension
$d=2$ with $\delta'=\delta''$, we obtain the following result for
interpoint distances in the unit disk.

\begin{theorem}
  \label{thr:2d}
  Assume that~(\ref{eq:su}) and (\ref{eq:ba}) hold with any fixed
  $\delta\in(0,\thf)$ uniformly over $u\in\SS^1$. Then
  (\ref{eq:sgamma}) holds. 
\end{theorem}

Instead of imposing (\ref{eq:su}) and (\ref{eq:ba}), it is possible to
deduce the limiting distribution in (\ref{eq:sgamma}) using a direct
assumption on the asymptotic distribution of $\Pit_{n\kappa}$.

\begin{theorem}
  \label{thr:altern}
  Assume that, for $[\delta',\delta'']\subset (0,\thf)$ with
  $\delta'<\delta''$,
  \begin{equation}
    \label{eq:new-con}
    \lim_{s\downarrow0}
    \frac{\Prob{\Pit_{n\kappa}\cap C_s(u,z_s)\neq \emptyset}}
    {\sigma(C(u,z_s))}=g(t,u)
  \end{equation}
  uniformly in $u\in\Sphere$ and $z_s\in [s^{\delta'},s^{\delta''}]$,
  where $n=ts^{-\gamma/2}$ and $\sigma$ s a finite measure on
  $\Sphere$. If the non-atomic part $\sigma'$ of $\sigma$ satisfies
  \textbf{(S)}, then
  \begin{multline}
    \label{eq:new-lim}
    \lim_{n\ti} \Prob{n^{2/\gamma}(2-\diam(\Pi_{n\kappa}))\leq t}\\
    =1-\exp\left\{-\thf \int_{\Sphere} g(t,u)\sigma'(du)\right\}
    \prod_{x_i\in\Sphere\atop \sigma(\{x_i\})>0}
    \Big(1-g(t,x_i)\sigma(\{x_i\})\Big)^{\thf} 
  \end{multline}
  for all $t\geq0$.
\end{theorem}
\begin{proof}
  For the proof we use the same sets partition and the sets $A^{(i)}$
  as in the proof of Theorem~\ref{thr:non-unif}. If $\sigma$ has an
  atomic part, choose the points $x_1,\dots,x_m$ in such a way that
  they have so many atoms of $\sigma$ among them that the total
  $\sigma$-content of the remaining atoms is less than $\eps$.
  
  In the remainder of the proof we need to split the product in the
  left-hand side of (\ref{eq:inf-product}) into the factors that
  correspond to the non-atomic and the atomic parts of $\sigma$.
  Notice that Lemma~\ref{le:large-block} is no longer needed to derive
  the asymptotics for $\Prob{\Pit(A_s^{(i)})=0}$ from the distribution
  of $\xit$. The square root of the product appears because we need to
  count only atoms from the half-sphere. Alternatively it is possible
  to take the produst only over $x_i\in\Sphere_+$.
\end{proof}

The cases when $\sigma$ has atoms often appear if $\kappa$ is the (say
uniform distribution) supported by a subset $K$ of $B$ and such that
$K$ is sufficiently ``sharply pointed'' near the points where its
diameter is achieved. The typical example of such $K$ is a segment,
see \ref{sec:distr-segm}. Other examples correspond to sets that
satisfy the conditions imposed in \cite{app:naj:rus02}.

\section{De-Poissonisation}
\label{sec:de-poissonisation}

Let $\Pi$ be the Poisson process with intensity measure $n\kappa$.
Given $\Pi(K)=n$, the distribution of $\Pi$ coincides with the
distribution of $\Xi_n=\{\xi_1,\dots,\xi_n\}$ being the binomial
process on $K$ that consists of i.i.d. points sampled from $\kappa$.
In the other direction, the distribution of $\Pi$ coincides with the
distribution of $\Xi_N$, where $N$ is the Poisson random variable of
mean $n$ independent of the i.i.d. points $\xi_i$'s distributed
according to $\kappa$.  This simple relationship makes it possible to
use the de-Poissonisation technique \cite{pen03} in order to obtain
the limit theorem for functionals of $P_n$ being the convex hull of
$\Xi_n$.  The key issue that simplifies our proofs is the monotonicity
of the diameter functional. Indeed, the diameter of $\Xi_n$ is
stochastically greater than the diameter of $\Xi_m$ for $n>m$.
Another useful tool is provided by the following lemma from
\cite[p.~18]{pen03}.

\begin{lemma}
  \label{le:poisson-tail}
  Let $N$ be a Poisson random variable with mean $\lambda$. For every
  $\gamma>0$ there exists a constant $\lambda_1=\lambda_1(\gamma)\geq
  0$ such that
  \begin{displaymath}
    \Prob{|N-\lambda|\geq \thf\lambda^{\thf+\gamma}} \leq
    2\exp\{-\frac{1}{9}\lambda^{2\gamma}\}
  \end{displaymath}
  for all $\lambda>\lambda_1$.
\end{lemma}

\begin{theorem}
  \label{th:depois}
  Let $\Psi:\sN\to\R$ be a monotonic functional defined on the space
  $\sN$ of finite subsets of $\R^d$. Furthermore, let $\Pi_{n \kappa}$
  be a Poisson process with intensity measure $n\kappa$ where $\kappa$
  is a probability measure on $\R^d$.  If, for some $\alpha>0$, the
  random variable $n^\alpha\Psi(\Pi_{n\kappa})$ converges in
  distribution to a random variable with cumulative distribution
  function $F$, then $n^\alpha\Psi(\Xi_n)$ also weakly converges to
  $F$, where $\Xi_n$ is a binomial process of $n$ i.i.d. points with
  common distribution $\kappa$.
\end{theorem}
\begin{proof}
  Without loss of generality assume that $\Psi$ is non-decreasing.
  Define $\gamma=\thf-\beta$ and $\eps_n=n^{-\beta}$ for some
  $\beta\in(0,\thf)$. By Lemma \ref{le:poisson-tail} and the
  monotonicity of $\Psi$,
  \begin{align*}
    \Prob{\Psi(\Pi_{n\kappa})\leq s} & \leq
    \Prob{\Psi(\Pi_{n\kappa})\leq s, \; |N-n|\leq n\eps_n}
    +\Prob{|N-n|>n\eps_n}\\
    & \leq \Prob{\Psi(\Xi_{n(1-\eps_n)})\leq s}
    +2\exp\{-\frac{1}{9}(2n)^{2\gamma}\}\,.
  \end{align*}
  for sufficiently large $n$.  Therefore, for every continuity point
  $t$ of $F$,
  \begin{align*}
    \lim_{n\ti}\Prob{\Psi(\Xi_n)n^\alpha \leq t}
    &= \lim_{n\ti}\Prob{\Psi(\xi_{n(1-\eps_n)})(n(1-\eps_n))^\alpha \leq t} \\
    &\geq \lim_{n\ti} \Prob{\Psi(\Xi_{n(1-\eps_n)})n^\alpha \leq t}\\
    &\geq \lim_{n\ti}\Prob{\Psi(\Pi_{n\kappa})n^\alpha \leq t}
    - 2 \exp\{-\frac{1}{9}(2n)^{2\gamma}\}\\
    &= \lim_{n\ti}\Prob{\Psi(\Pi_{n\kappa})n^\alpha \leq t}\\
    &= F(t)\,.
  \end{align*}
  Similarly, we have
  \begin{align*}
    \Prob{\Psi(\Pi)\leq s}&\geq
    \Prob{\Psi(\Xi_{n(1+\eps_n)})\leq s}\Prob{|N-n|\leq n\eps_n}\\
    &\geq \Prob{\Psi(\Xi_{n(1+\eps_n)})\leq s} -
    2\exp\{-\frac{1}{9}(2n)^{-2\gamma}\}\,,
  \end{align*}
  so that
  \begin{align*}
    \lim_{n\ti}\Prob{\Psi(\Xi_n)n^\alpha \leq t}
    &= \lim_{n\ti}\Prob{\Psi(\xi_{n(1+\eps_n)})(n(1+\eps_n))^\alpha
      \leq t} \\
    &\leq \lim_{n\ti} \Prob{\Psi(\Xi_{n(1+\eps_n)})n^\alpha \leq t}\\
    &\leq \lim_{n\ti}\Prob{\Psi(\Pi_{n\kappa})n^\alpha \leq t}
    - 2 \exp\{-\frac{1}{9}(2n)^{2\gamma}\}\\
    &= \lim_{n\ti}\Prob{\Psi(\Pi_{n\kappa})n^\alpha \leq t}\\
    &= F(t).
  \end{align*}
\end{proof}

In particular, Theorem~\ref{th:depois} is applicable for the
functional $\Psi(\Xi_n)=2-\diam(\Xi_n)$, so that all results
available for diameters of Poisson processes can be immediately
reformulated for binomial processes.

\section{Spherically symmetric distributions}
\label{sec:spherically-symmetric}

Let $\xi_1,\dots,\xi_n$ be independent points distributed according to
a spherically symmetric (also called ``isotropic'') density $\kappa$
restricted on $B$.  Spherically symmetric distributions are closed
with respect to convolution, so that $\xit=\xi_1-\xi_2$ is spherically
symmetric too.  Therefore, $\|\xit\|$ and $\xit/\|\xit\|$ are
independent, see e.g.  \cite{gne98}.  Then, for any measurable
$A\subset\Sphere$
\begin{displaymath}
  \Prob{\xit \in A_s}
  = \Prob{\|\xit\| \ge 2-s}
  \frac{\mu_{d-1}(A)}{\mu_{d-1}(\Sphere)}\,.
\end{displaymath}
Therefore (\ref{eq:su}) is fulfilled if, for some $\gamma>0$,
\begin{equation}
  \label{eq:su-spherically-symmetric}
  \lim_{s \to 0} \Prob{\|\xit\| \ge 2-s}s^{-\gamma}
  = \sigma_0\in(0,\infty)\,,
\end{equation}
where the limit $\sigma_0$ then becomes the total mass of $\sigma$, so
that $\sigma$ is the surface area measure on $\Sphere$ normalised to
have the total mass $\sigma_0$.

Furthermore, (\ref{eq:ba}) holds if
\begin{equation}
  \label{eq:ba-spherically-symmetric}
  \lim_{s\downarrow 0}s^{\delta(d-1)-\gamma/2}\Prob{\eta\leq s}=0\,,
\end{equation}
where $\eta=1-\|\xi_1\|$.  

\begin{lemma}
  \label{lem:zeta}
  If $\eta_1$ and $\eta_2$ are i.i.d. distributed as $1-\|\xi_1\|$ and
  $\zeta=\eta_1+\eta_2$, then
  \begin{equation}
    \label{eq:zeta-tail}
    \lim_{s\downarrow0} \frac{\Prob{\|\xit\|\geq 2-s}}
    {\E((s-\zeta)^{(d-1)/2}\Ind_{\zeta\leq s})}
    = \frac{2^{d-1}\Gamma(\frac{d}{2})}{(d-1)\pi^\thf\Gamma(\frac{d-1}{2})}\,.
  \end{equation}
\end{lemma}
\begin{proof}
  By the cosine theorem and the fact that $\xit$ has the same
  distribution as $\xi_1+\xi_2$, we write
  \begin{displaymath}
    \Prob{\|\xit\| \geq 2 -s} = \Prob{\|\xi_1\|^2+\|\xi_2\|^2 +
      2\|\xi_1\|\|\xi_2\|\cos\beta \geq (2 -s)^2},
  \end{displaymath}
  where $\beta$ denotes the angle between $\xi_1$ and $\xi_2$.  Hence,
  \begin{displaymath}
    \Prob{\|\xit\| \geq 2-s} = \Prob{\cos\beta \geq 1 - q}\,,
  \end{displaymath}
  where
  \begin{displaymath}
    q = \frac{(2-\zeta)^2-(2-s)^2}{2 \|\xi_1\|\|\xi_2\|}\,.
  \end{displaymath}
  If $q\geq 0$ (i.e. $\zeta\leq s$)
  \begin{displaymath}
    \Prob{\|\xit\| \geq 2 -s}
    = \thf\Prob{\cos^2\beta \geq(1-q)^2, \zeta \leq s}
    = \thf \E\Big(\int_{(1-q)^2}^1 f(t)dt\Ind_{\zeta \leq s}\Big)\,,
  \end{displaymath}
  where the probability density function
  \begin{displaymath}
    f(t) = \frac{\Gamma(\frac{d}{2})}{\pi^\thf\Gamma(\frac{d-1}{2})} 
    t^{-\thf}(1-t)^{\frac{d-3}{2}}\,, \quad t\in[0,1]\,,
  \end{displaymath}
  of $\cos^2\beta$ corresponds to the Beta-distribution with
  parameters $\thf$ and $(d-1)/2$, see \cite[Prop.~2]{mat:ruk93}.
  Substituting $x = 1 - t$ leads to
  \begin{displaymath}
    \Prob{\|\tilde\xi\| \geq 2 -s} = c_1 
    \E\Big(\int_0^{2q(1-\frac{q}{2})}(1-x)^{-\thf}
    x^{\frac{d-1}{2}-1}dx\Ind_{\zeta \leq s}\Big)\,,
  \end{displaymath}
  where
  \begin{displaymath}
    c_1 = \thf
    \frac{\Gamma(\frac{d}{2})}{\pi^\thf\Gamma(\frac{d-1}{2})}\,.
  \end{displaymath}
  The inequality
  \begin{displaymath}
    1 \leq (1-x)^{-\thf} \leq (1 - q)^{-1} \leq
    ((1-s)^2-2s)^{-1}
  \end{displaymath}
  leads to the bounds
  \begin{equation}
    \label{eq:p}
    c_1\E(I \Ind_{\zeta \leq s}) \leq \Prob{\|\xit\| \geq 2 -s} \leq
    c_1((1-s)^2-2s)^{-1}\E(I\Ind_{\zeta \leq s})\,,
  \end{equation}
  where 
  \begin{displaymath}
    I = \int_0^{2q(1-\frac{q}{2})}x^{\frac{d-1}{2}-1}dx = 
    \frac{2}{d-1}(2q)^\frac{d-1}{2}(1-\frac{q}{2})^\frac{d-1}{2}.
  \end{displaymath}
  By the fact that
  \begin{displaymath}
    1 \leq (\|\xi_1\|\|\xi_2\|)^{-1}\leq (1-s)^{-2}
  \end{displaymath}
  and
  \begin{displaymath}
    1 - \frac{s}{(1-s)^2} \leq
    1-\frac{q}{2} \leq 1\,,
  \end{displaymath}
  we further get the bounds
  \begin{multline*}
    \frac{2}{d-1}((2-\zeta)^2-(2-s)^2)^\frac{d-1}{2} \left(1 - 
      \frac{s}{(1-s)^2}\right)^\frac{d-1}{2}\leq I \\
    \leq \frac{2}{d-1}((2-\zeta)^2-(2-s)^2)^\frac{d-1}{2}
    (1-s)^{-(d-1)}\,.
  \end{multline*}
  Since
  \begin{displaymath}
    (4-2s)(s-\zeta) \leq (2 - \zeta)^2-(2-s)^2 \leq 4(s-\zeta)\,,
  \end{displaymath}
  the following bounds for $I$ hold
  \begin{multline*}
    \frac{2^d}{d-1}(s-\zeta)^\frac{d-1}{2} \left(1 - 
      \frac{s}{(1-s)^2}\right)^\frac{d-1}{2}
    (1-\frac{s}{2})^\frac{d-1}{2}\leq I \\
    \leq \frac{2^d}{d-1}(s-\zeta)^\frac{d-1}{2}(1-s)^{-(d-1)}\,.
  \end{multline*}
  Plugging these bounds in (\ref{eq:p}) yields the result.
\end{proof}

The following result settles the case when the density of $\eta$ is
equivalent to a power function for small arguments.

\begin{theorem}
  \label{thr:spherical-gen}
  Assume that $d\geq2$ and for some $\alpha\geq0$ the cumulative
  distribution function $F(x)=\Prob{\eta\leq x}$ of $\eta$ satisfies
  \begin{equation}
    \label{eq:s-eta}
    \lim_{s\downarrow0} s^{-\alpha}F(s) = a\in(0,\infty)\,.
  \end{equation}
  Then 
  \begin{equation}
    \label{eq:s-limit}
    \lim_{n\ti} \Prob{n^{2/\gamma}(2-\diam(\Pi_{n\kappa}))\leq t}
    =1-e^{-\thf t^\gamma\sigma_0}\,, \quad t\geq0\,,
  \end{equation}
  where $\gamma=\thf(d-1)+2\alpha$ and 
  \begin{displaymath}
    \sigma_0=a^2c
    \begin{cases}
      1\,, & F(0)>0\,,\\
      \alpha^2 \Gamma(\alpha)^2
      \frac{\Gamma(\thf(d+1))}{\Gamma(2\alpha+\thf(d+1))}\,,
      & F(0)=0\,,
    \end{cases}
  \end{displaymath}
  with $c$ given by the right-hand side of (\ref{eq:zeta-tail}). 
\end{theorem}
\begin{proof}
  The integration by parts leads to
  \begin{multline*}
    \E((s-\zeta)^{(d-1)/2}\Ind_{\zeta\leq s})
    =F(0)^2s^{(d-1)/2}\\
    +\frac{(d-1)(d-3)}{4}
    \int_0^s\int_0^{s-x_1} F(x_1)F(x_2)(s-x_1-x_2)^{(d-5)/2}dx_1dx_2\,.
  \end{multline*}
  If $F(0)>0$, then (\ref{eq:s-eta}) implies that $\alpha=0$, so that
  (\ref{eq:su-spherically-symmetric}) holds with $\gamma=\thf(d-1)$
  and $\sigma_0=F(0)^2c=a^2c$ by Lemma~\ref{lem:zeta}.
  
  If $F(0)=0$, then (\ref{eq:s-eta}) yields that
  $\E((s-\zeta)^{(d-1)/2}\Ind_{\zeta\leq s})$ is equivalent as $s\downarrow0$
  to
  \begin{align*}
    s^\gamma a^2 \frac{(d-1)(d-3)}{4}
    &\int_0^1\int_0^{1-t_1} t_1^\alpha
    t_2^\alpha(1-t_1-t_2)^{(d-5)/2}dt_1dt_2 \\
    &=s^\gamma a^2 \frac{(d-1)(d-3)}{4}
    \Beta({\small\alpha+1,\alpha+\frac{d-1}{2}})
    \Beta({\small\alpha+1,\frac{d-3}{2}})\\   
    &=s^\gamma a^2 \alpha^2 \Gamma(\alpha)^2
    \frac{\Gamma(\thf(d+1))}{\Gamma(2\alpha+\thf(d+1))}
  \end{align*}
  with $\gamma=\thf(d-1)+2\alpha$.  Finally,
  (\ref{eq:su-spherically-symmetric}) follows from
  Lemma~\ref{lem:zeta}. It remains to show that
  (\ref{eq:ba-spherically-symmetric}) holds, i.e.
  \begin{displaymath}
    \delta(d-1)-\thf\gamma+\alpha >0\,.
  \end{displaymath}
  Using the expression for $\gamma$, it suffices to note that
  $\delta(d-1)-\frac{1}{4}(d-1)>0$ if $\delta\in(\frac{1}{4},\thf)$,
  so it is possible to choose $[\delta',\delta'']\subset
  (\frac{1}{4},\thf)$.
\end{proof}

It should be noted that (\ref{eq:s-eta}) can be replaced by the
requirement that $F$ is regular varying at zero. However, in this case
the constants involved in the formula for $\sigma_0$ are given by the
integrals of the slowly varying part of $F$.

Using similar arguments, it is possible to check (\ref{eq:su}) and
(\ref{eq:ba}) if $\xi=\eta u$ for independent $\eta$ and $u$, where
$\eta$ distributed on $[0,1]$ and $u$ is distributed on $\Sphere$.

\section{Examples}
\label{sec:examples}

\subsection{Uniformly distributed points in the ball}

Consider the case of random points uniformly distributed in $B$.

\begin{theorem}
  \label{theorempoisson}
  As $n\to\infty$, the diameter of the convex hull of a homogeneous
  Poisson process $\Pi_\lambda$ with intensity $\lambda=n/\mu_d(B)$
  restricted on the $d$-dimensional unit ball, $d \ge 2$, has limit
  distribution
  \begin{displaymath}
    \Prob{n^{\frac{4}{d+3}}(2 - \diam \Pi_\lambda)\le t}
    \to 1 - \exp\left\{-\thf\, c\, t^{\frac{d+3}{2}}\right\}, \quad t > 0\,,
  \end{displaymath}
  where
  \begin{equation}
    \label{eq:c-unif}
    c= \frac{2^{d+1} d \Gamma(\frac{d}{2}+1)}
      {\sqrt{\pi}(d+1)(d+3)\Gamma(\frac{d+1}{2})}\,.
  \end{equation}
\end{theorem}
\begin{proof}
  The tail behaviour of $\|\xi_1\|$ is determined by
  \begin{displaymath}
    \Prob{\|\xi_1\| \geq 1 - s} = 1 - \frac{\mu_d(B_{1-s})}{\mu_d(B)}
    = 1 - (1-s)^d\,,
  \end{displaymath}
  so that Theorem~\ref{thr:spherical-gen} is applicable with $\alpha=1$
  and $a=d$. 
\end{proof}

By the de-Poissonisation argument, Theorem~\ref{theorempoisson} yields
Theorem~\ref{theoremuniform}.  Note that in case $d=2$ the constant
$c$ equals $16/(15\pi)$, which also corresponds to the bounds given in
(\ref{eq:ap-r}). The tail behaviour of $\|\xit\|$ can also be obtained
from the explicit formula for the distribution of the length of a
random chord in the unit ball, see \cite[2.48]{ke:mor}.

\subsection{Uniformly distributed points on the sphere}

Another example of a spherically symmetric distribution is given by
the uniform distribution on $\Sphere$, i.e. if 
\begin{displaymath}
  \kappa(A) = \mu_{d-1}(A)/\mu_{d-1}(\Sphere)
\end{displaymath}
for all measurable $A\subset\Sphere$. The following result follows
from Theorem~\ref{thr:spherical-gen} in case $a=F(0)=1$ and
$\alpha=0$. 

\begin{theorem}
  \label{thr:boundary}
  If $\Pi$ is the homogeneous Poisson process on $\Sphere$ with the
  total intensity $n$, then for any $d \geq 2$
  \begin{displaymath}
    \lim_{n\ti} \Prob{n^\frac{4}{d-1}(2-\diam(\Pi))\leq t}
    =1-\exp\left\{-\thf\, c\, t^\frac{d-1}{2}\right\}\,,\quad t>0\,,
  \end{displaymath}
  where
  \begin{equation}
    \label{eq:c-sphere}
    c=\frac{2^{d-1}\Gamma(\frac{d}{2})}{(d-1)\pi^\thf\Gamma(\frac{d-1}{2})}\,.
  \end{equation}
\end{theorem}

Alternatively, the tail behaviour of $\|\xit\|$ may be derived from
the explicit formula for the distribution of the distance between two
uniform points on the unit sphere, see \cite{ala76}. 

Similarly, it is possible to obtain limit theorems for a spherically
symmetric $\xi$ in case the norm $\|\xi\|$ has a rather general
distribution which is regular varying near its right end-point being
1.

\subsection{Distribution in spherical sectors}
\label{sec:prop-distr}

This section provides a simple example, where $\kappa$ is not
spherically symmetric. Consider some spherically symmetric measure
$\kappa'$ which satisfies (\ref{eq:su-spherically-symmetric}) and
(\ref{eq:ba-spherically-symmetric}) for some $c$ and $\gamma$, and a
spherical sector $L$ defined by
\begin{displaymath}
  L = \{tx: x \in A, t \in [-1, 1]\}
\end{displaymath}
for some fixed $r>0$, where $A$ is a spherically convex subset of the
unit sphere. If $L\neq B$, then 
\begin{displaymath}
  \kappa(S) = \kappa'(S \cap L)/\kappa'(L)
\end{displaymath}
defines a not spherically symmetric measure for all measurable $S
\subset B$. Denote by $\xi_1$ and $\xi_2$ two independent points
sampled from $\kappa$ and by $\xi'_1, \xi'_2$ two independent points
distributed according to $\kappa'$, respectively. By the construction
of $\kappa$ and the spherical symmetry of $\kappa'$ we can write
\begin{equation}
  \label{eq:propeller}
  \Prob{\xit\in C_s(u,r)} = \Prob{\|\xi'_1-\xi'_2\| \geq
    2-s}\frac{\mu_{d-1}(C(u,r))\cap A)}{\mu_{d-1}(A)}
\end{equation}
for all spherical balls $C_s(u,r)$. For all measurable $F\subset
\Sphere$ define
\begin{displaymath}
  \sigma(F) = c\, \frac{\mu_{d-1}(F\cap A)}{\mu_{d-1}(\Sphere)}\,.
\end{displaymath}
It is easy to verify that $\sigma$ satisfies condition (\ref{eq:su}),
since by (\ref{eq:propeller})
\begin{displaymath}
  \lim_{s\downarrow0} \sup_{u \in \supp \sigma}
  \left|\frac{\Prob{\xit\in C_s(u,r)}}{s^\gamma\sigma(C(u,r))}
    -1\right|=0\,.
\end{displaymath}
Since
\begin{displaymath}
  \kappa(D_s(C_{s}(u,s^\delta))) =
  \kappa'(D_s(C_{s}(u,s^\delta))\cap L)/\kappa'(L)\,,
\end{displaymath}
condition (\ref{eq:ba}) is fulfilled and Theorem~\ref{thr:non-unif}
holds with
\begin{displaymath}
  \sigma_0=c\frac{\mu_{d-1}(A)}{\mu_{d-1}(\Sphere)}\,.
\end{displaymath}

\subsection{Non-uniform angular distributions}
\label{sec:non-uniform-angular}

Assume that $\xi$ is distributed on the boundary of the unit circle in
$\R^2$ according to some not necessarily symmetrical probability
measure $\kappa$, which can be then considered a measure on
$[0,2\pi)$. If $\xi_1$ and $\xi_2$ are distributed on $[0,2\pi)$
according to $\kappa$, then
\begin{align*}
  \Prob{1-\cos(\xi_1-\xi_2)\leq 2&s(1-s/2)(1-s)^{-2},\,
    |\xi_1+\xi_2-2u|\leq 2s^\delta}\\
  &\leq \Prob{\xit\in C_s(u,s^\delta)}\\
  &\leq \Prob{1-\cos(\xi_1-\xi_2)\leq 2s,\,
    |\xi_1+\xi_2+\pi-2u|\leq 2s^\delta}\,,
\end{align*}
where the addition of angles is understood by modulus $2\pi$.
Thus, $\Prob{\xit\in C_s(u,s^\delta)}$ is equivalent as $s\downarrow 0$ to 
\begin{displaymath}
  \Prob{|\xi_1+\xi_2|\leq 2\sqrt{s},\, |\xi_1+\xi_2+\pi-2u|\leq 2s^\delta}\,.
\end{displaymath}
Assume that the distribution $\kappa$ has bounded density $f$ with
respect to the length measure on the unit circle. Then the probability
above is equivalent to
\begin{displaymath}
  2\sqrt{s}4s^\delta f(u)f(u+\pi)=4s^{\gamma} (2s^\delta)f(u)f(u+\pi)\,.
\end{displaymath}
Thus, (\ref{eq:su}) holds with $\gamma=\thf$ and
\begin{displaymath}
  \sigma_0=4\int_0^{2\pi} f(u)f(u+\pi)du\,.
\end{displaymath}
The boundedness of $f$ also implies that (\ref{eq:ba}) holds, so that
the limit distribution is given by (\ref{eq:sgamma}).  In particular
if $\kappa$ is uniform on the circle, then $f(u)=1/(2\pi)$, so that
$\sigma_0=4/(2\pi)=2/\pi$, so that
\begin{displaymath}
  \lim_{n\ti} \Prob{n^4(2-\diam(\Pi_{n\kappa}))\leq t}
    =1-e^{-\frac{1}{\pi}\sqrt{t}}\,,\quad t\geq0\,,
\end{displaymath}
which also corresponds to the result of Theorem~\ref{thr:boundary} for
$d=2$.

\subsection{Segments and disks in the unit ball}
\label{sec:distr-segm}

Assume that $L_1,\dots,L_m$ are segments that obtained by intersection
the unit ball with $m$ different lines passing through the origin.
Assume that $\Prob{\xi\in L_i}=p_i$, $i=1,\dots,m$, and given $\xi\in
L_i$, $\xi$ is distributed according to the length measure on $L_i$.
If $L_i=[-x_i,x_i]$, then let $\sigma$ be an atomic measure with unit
atoms at $\{\pm x_i, \,i=1,\dots,m\}$. The one-dimensional result for
the range of a uniform random sample \cite[???]{gal78} implies that
\begin{displaymath}
  \lim_{s\downarrow 0}
  \frac{\Prob{\Pit_{n\kappa}\cap C_s(x_i,z_s)\neq \emptyset}}
    {\sigma(C(x_i,z_s))}=1-\Big(1+\thf tp_i\Big)e^{-\thf tp_i}\,,
\end{displaymath}
where $n=t/s$, i.e. $\gamma=2$. Theorem~\ref{thr:altern} implies that 
\begin{equation}
  \label{eq:point-seg}
  \lim_{n\ti} \Prob{n(2-\diam(\Pi_{n\kappa}))\leq t}
  =1-e^{-\thf t} \prod_{i=1}^m \Big(1+\thf tp_i\Big)\,.
\end{equation}

If $L_1,\dots,L_m$ are obtained as intersections of the unit ball with
linear subspaces of possibly different dimensions, then only those
with the smallest dimension of these subspaces contribute to the
asymptotic distribution of the maximum interpoint distance. If the
smallest dimension is at least $2$, then $\sigma$ is non-atomic and
Theorem~\ref{thr:non-unif} is applicable as in the case of uniformly
distributed points. Otherwise, we arrive at the formula above for
segments.

\query{more here}

Assume now that the number of atoms $L_i=[-x_i,x_i]$, $i\geq1$.
Without loss of generality assume that $x_i\to x_0$ as $i\ti$ and
$x_i\neq x_0$ for all $i$.  Let $\nu$ be the measure on $\Sphere$ with
atoms $\pm x_i$ with $\nu(\{x_i\})=\nu(\{-x_i\})=p_i$, $i\geq1$. In
comparison with the case of a finite number of segments, we need also
to find the limit of $\Prob{\Pit_{n\kappa}\cap C_s(x_0,z_s)\neq
  \emptyset}$ as $s\to0$. Notice that $\nu(C_s(x_0,z_s))=q_s\to0$ as
$s\downarrow0$. Thus, $\Prob{\Pit_{n\kappa}\cap C_s(x_0,z_s)\neq
  \emptyset}$ is bounded above by the probability that the Poisson
point process with the total intensity $nq_s$ on $[-x_0,x_0]$ has
diameter that exceeds $2-s$. Using one-dimensional result, it is easy
to see that the corresponding limit is zero if $\gamma=2$. In order to
arrive at a non-trivial limit, we need to set $\gamma>2$, which is
impossible, since the normalisation $n^\gamma$ is to big for the
diameters of the Poisson processes restricted on the individual
segments $L_i$, $i\geq1$. Therefore, (\ref{eq:point-seg}) holds in
this case with the infinite product, i.e. for $m=\infty$. 

Foe instance, assume that $p_i=\zeta(2)i^{-2}$, $i\geq1$, where
$\zeta$ is the zeta-function. Using a formula for infinite product 
\cite[(89.5.16)]{han75} we obtain
\begin{displaymath}
  \lim_{n\ti} \Prob{n(2-\diam(\Pi_{n\kappa}))\leq t}
  =1-e^{-\thf t} \frac{1}{\pi}\sqrt{\frac{2\zeta(2)}{t}}\; 
  \sinh \pi\sqrt{\frac{t}{2\zeta(2)}}\,.
\end{displaymath}

\newcommand{\noopsort}[1]{} \newcommand{\printfirst}[2]{#1}
  \newcommand{\singleletter}[1]{#1} \newcommand{\switchargs}[2]{#2#1}

\end{document}